\documentclass[11pt, a4paper]{amsart}

\usepackage[T1]{fontenc}
\usepackage[utf8]{inputenc}
\usepackage{lmodern}

\usepackage[american]{babel}
\usepackage[babel, final]{microtype}
\usepackage{amssymb}

\usepackage{mathtools}

\usepackage{ifpdf}
\usepackage{comment}
\usepackage{multirow}

\usepackage[pdftex, dvipsnames]{xcolor}
\usepackage[pdftex, final]{graphicx}
\usepackage{caption}
\usepackage{subcaption}
\usepackage{pinlabel}
\usepackage[pdftex,%
    a4paper,%
    includehead,%
    includefoot,%
    nomarginpar,%
    lmargin=1in,%
    rmargin=1in,%
    tmargin=1in,%
    bmargin=1in,%
]{geometry}
\usepackage[pdftex,%
    final,%
    colorlinks=true,%
    linkcolor=NavyBlue,%
    citecolor=NavyBlue,%
    filecolor=NavyBlue,%
    menucolor=NavyBlue,%
    urlcolor=NavyBlue,%
    bookmarks=true,%
    bookmarksdepth=3,%
    bookmarksnumbered=true,%
    bookmarksopen=true,%
    bookmarksopenlevel=2,%
]{hyperref}
\hypersetup{
    pdftitle={Deep and shallow slice knots in 4-manifolds},
    pdfauthor={Michael R. Klug, Benjamin M. Ruppik}
}

\usepackage{overpic}

\usepackage{tikz}
\usetikzlibrary{tikzmark}

\usepackage{todonotes}
\usepackage{marginnote}

\input{macros}

\newcommand{\abs}[1]{\mathchoice%
  {\left\lvert #1 \right\rvert}%
  {\lvert #1 \rvert}%
  {\lvert #1 \rvert}%
  {\lvert #1 \rvert}%
}

\DeclareMathOperator{\interior}{int}
\DeclareMathOperator{\sign}{sign}

\newcommand{\sphere}[1]{{S}^{#1}}
\newcommand{\disk}[1]{{D}^{#1}}
\newcommand{\ball}[1]{{B}^{#1}}
\newcommand{\interval}{{I}}

\newcommand{\CP}[1]{\mathbb{CP}^{#1}}

\newcommand{\laurent}{\ZZ[t, t^{-1}]}

\newcommand{\ZZ}{\mathbb{Z}}

\newcommand{\SPC}{SPC4 }

\title{
Deep and shallow slice knots in 4-manifolds
}

\author{Michael R. Klug}
\address{UC Berkeley math department}
\email{\href{mailto:michael.r.klug@gmail.com}{michael.r.klug@gmail.com}}
\urladdr{\url{https://math.berkeley.edu/~mrklug/}}

\author{Benjamin M. Ruppik}
\address{Max-Planck-Institut f{\"ur} Mathematik, Bonn, Germany}
\thanks{MK and BR were supported by the Max Planck Institute for Mathematics in Bonn.}
\email{\href{mailto:bruppik@mpim-bonn.mpg.de}{bruppik@mpim-bonn.mpg.de}}
\urladdr{\url{https://ben300694.github.io/}}

\def\subjclassname{\textup{2020} Mathematics Subject Classification}
\expandafter\let\csname subjclassname@1991\endcsname=\subjclassname
\expandafter\let\csname subjclassname@2000\endcsname=\subjclassname
\subjclass{57K40 
    (Primary),
    57K10 
    (Secondary);
    \hfill
    Date: \today
}


\begin{document}
 
\maketitle

\begin{abstract}
    We consider slice disks for knots
    in the boundary of a smooth compact 4-manifold $X^{4}$.
    We call a knot $K \subset \partial X$ \textit{deep slice} in $X$
    if there is a smooth properly embedded $2$-disk in $X$
    with boundary $K$, but $K$ is not concordant
    to the unknot in a collar neighborhood
    $\partial X \times \interval$
    of the boundary.
    
    We point out how this concept relates to various well-known conjectures
    and give some criteria for the nonexistence of such deep slice knots.
    Then we show, using the Wall self-intersection invariant
    and a result of Rohlin,
    that every 4-manifold consisting of
    just one 0- and a nonzero number of
    2-handles always has a deep slice knot in the boundary.
    
    We end by considering 4-manifolds where every knot
    in the boundary bounds an embedded disk in the interior.
    A generalization of the Murasugi-Tristram inequality 
    is used to show that there does not exist
    a compact, oriented $4$-manifold $V$
    with spherical boundary
    such that every knot $K \subset \sphere{3} = \partial V$
    is slice in $V$ via a null-homologous disk.
\end{abstract}

\section{Overview}

The Smooth 4-Dimensional Poincar\'e Conjecture (SPC4)
proposes that every
closed smooth 4-manifold $\Sigma$ that is homotopy equivalent to $\sphere{4}$
is diffeomorphic to the standard $\sphere{4}$.
By work of Freedman \cite{freedman1982topology},
it is known that if $\Sigma$ is homotopy equivalent to $\sphere{4}$,
then $\Sigma$ is in fact homeomorphic to $\sphere{4}$.
In stark contrast to the SPC4,
it might be the case that every compact smooth 4-manifold
admits infinitely many distinct smooth structures.
The existence of an exotic homotopy 4-sphere is equivalent
to the existence of an exotic
contractible compact manifold with
$\sphere{3}$ boundary \cite[p.\ 113]{MR0190942},
henceforth called an \emph{exotic homotopy 4-ball}.

One possible approach to proving that a proposed exotic homotopy 4-ball $\mathcal{B}$
is in fact exotic is to find a knot $K \subset \sphere{3} = \partial \mathcal{B}$,
such that there is a smooth properly embedded disk
$\disk{2} \hookrightarrow \mathcal{B}$,
with $\partial \disk{2}$ mapped to $K$,
where $K$ is not smoothly slice in the usual sense
in the standard 4-ball $\ball{4}$.
A knot is (topologically/smoothly) slice in $\ball{4}$ 
if and only if it is null-concordant in
$\sphere{3} \times \interval = \sphere{3} \times [0, 1]$,
i.e.\ there is a properly embedded (locally flat/smooth) cylinder
$\sphere{1} \times \interval \hookrightarrow \sphere{3} \times \interval$
whose oriented boundary is
$K \subset \sphere{3} \times \{ 0 \}$ together with
the unknot $U \subset \sphere{3} \times \{ 1 \}$.
Another way of thinking about this strategy is
that we want to find a knot $K$ in
$\sphere{3} = \partial \mathcal{B}$
that bounds a properly embedded smooth disk in $\mathcal{B}$
but does not bound any such disk that is contained in a 
collar $\sphere{3} \times \interval$ of the boundary of $\mathcal{B}$.
In this case, to verify the sliceness of $K$, we have to go ``deep'' into $\mathcal{B}$.  

An easier task might be to find a homology 4-ball $X$
with $\sphere{3}$ boundary such that there is a knot in the boundary
that bounds a smooth properly embedded disk in $X$ but not in $\ball{4}$,
however, this is also an open problem.
In \cite{Freedman_2010}, the authors investigate the possibility of proving
that a homotopy 4-ball $\mathcal{B}$ with $\sphere{3}$ boundary is exotic by taking a
knot in the boundary that bounds a smooth properly embedded disk in $\mathcal{B}$
and computing the $s$-invariant of $K$, in the hopes that $s(K) \neq 0$,
whereby they could then conclude that $\mathcal{B}$ is exotic.
Unfortunately for this approach as noted in the paper,
it turns out that the homotopy 4-ball that they were studying
was in fact diffeomorphic to $\ball{4}$, see \cite{akbulut2010cappell}.
It is still open whether the $s$-invariant can obstruct the
sliceness of knots in $\ball{4}$ that are slice in some homotopy 4-ball,
as is noted in the corrigendum to \cite{kronheimer2013gauge}.

Motivated by this, we make the following definitions:
For a 3-manifold $M^3$ containing a knot
$K \colon \sphere{1} \hookrightarrow M$, we say
that \emph{$K$ is null-concordant in $M \times \interval$}
if there is a smoothly properly embedded annulus
$\sphere{1} \times \interval \hookrightarrow M \times I$
cobounding $K \subset M \times \{ 0 \}$ on one end
and an unknot contained in a 3-ball $U \subset B^3 \subset M \times \{ 1 \}$
on the other.
Equivalently, $K \subset M \times \{ 0 \}$ bounds a smoothly properly
embedded disk in $M \times \interval$.

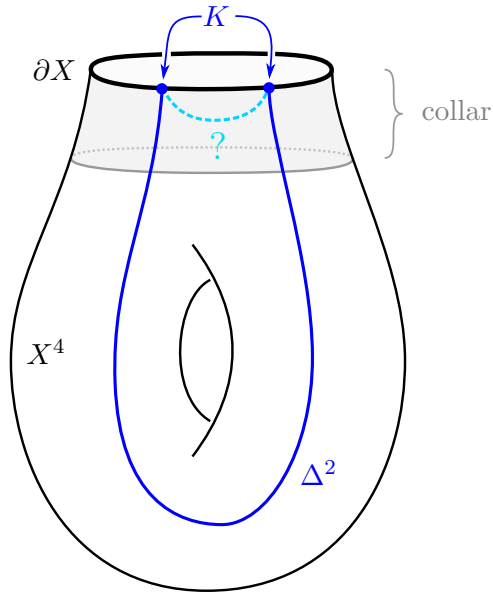
\begin{figure}
    \centering
    \begin{overpic}[width=0.35\textwidth]{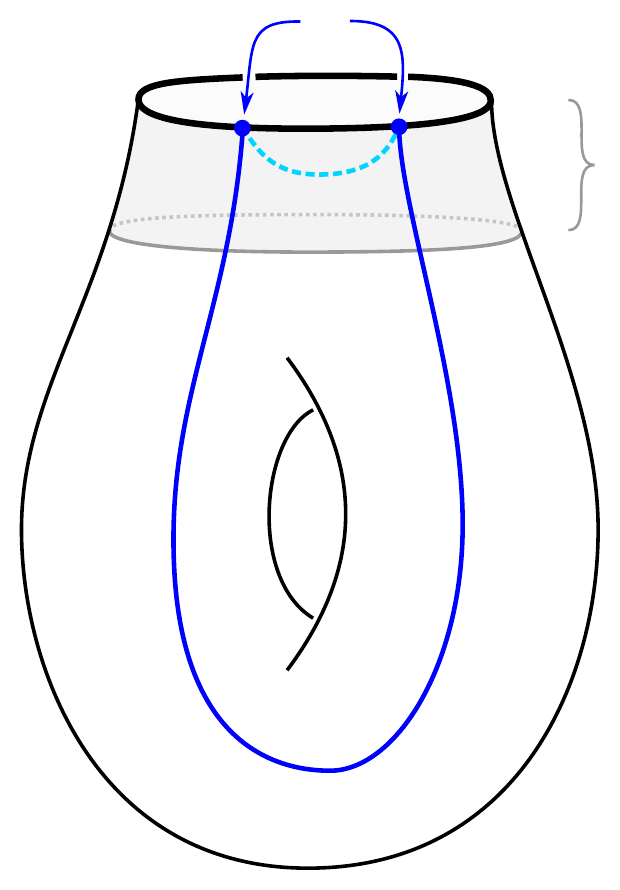}
        \definecolor{blue-green}{rgb}{0.0, 0.87, 0.87}
        \put(34, 96){\color{blue} $K$}
        \put(6, 87){$\partial X$}
        \put(35, 74){\Large \color{blue-green} ?}
        \put(70, 80.5){\color{gray} collar}
        \put(5, 40){$X^4$}
        \put(50, 20){\color{blue} $\Delta^2$}
    \end{overpic}
    \caption{
        Schematic of a deep slice disk $\Delta^2$ (blue)
        in a 4-manifold $X^4$,
        with boundary the knot $K \subset \partial X$.
        The knot $K$ is called \emph{deep slice}
        if it does not bound a properly embedded disk in a collar
        of the boundary, indicated by the (light blue) dashed lines.
    }
    \label{fig:deep_shallow_schematic}
\end{figure}

\begin{definition}[Deep slice/Shallow slice]
    Let $X^{4}$ be a smooth compact 4-manifold with nonempty boundary $\partial X$.
    We call a knot $K \subset \partial X$ \emph{slice} in $X$ if there is a smooth properly embedded disk in $X$
    with boundary $K$.
    We call a knot $K \subset \partial X$ \textit{shallow slice} in $X$
    if there is a smooth properly embedded disk in $\partial X \times I$
    with boundary $K$ -- this is equivalent
    to $K$ being null-concordant in the collar $\partial X \times \interval$.
    If $K$ is slice in $X$ but not shallow slice,
    we will call it \textit{deep slice} in $X$.
    See \autoref{fig:deep_shallow_schematic} for a schematic
    illustration of these definitions.
\end{definition}

In this language,
Problem 1.95 on Kirby's list \cite{kirby1995problems}
(attributed to Akbulut)
can be reformulated as follows:
Are there contractible smooth 4-manifolds
with boundary an integral homology 3-sphere
which contain deep slice knots that are null-homotopic in the boundary?
Note that any knot that is not nullhomotopic in the boundary
will not be shallow slice and
thus if it is slice, it will be deep slice.
For this reason we will be looking for deep slice knots
that are null-homotopic in the boundary.
We will often consider our knots
to be contained in 3-balls in the boundary,
which we call \emph{local knots},
so we can freely consider them in the boundary of any 4-manifold and
discuss if they are slice there.
To avoid confusion
when we say that a (local) knot in a 3-manifold $M^3$ is slice
we will usually qualify it with ``in $X^4$''.

\subsection{Outline}

In the first part of this paper
we will restrict ourselves to the smooth category,
starting in \autoref{sec:nonexistence_deep_slice}, 
where we discuss a condition that guarantees that some 4-manifolds
have no deep slice knots and related results.
In \autoref{sec:existence_deep_slice}, we prove that every 2-handlebody
has a deep slice knot in its boundary.
To do this we employ the Wall self-intersection number
and a result of Rohlin
which we discuss briefly.

In \autoref{sec:universal_slicing}, we recall the Norman-Suzuki trick
and observe that every 3-manifold bounds a 4-manifold
where every knot in the boundary bounds a properly embedded disk.
In contrast, if we restrict to slice disks trivial in
relative second homology, 
we will see that every compact topological $4$-manifold
with boundary $\sphere{3}$
contains a knot which
does not bound a null-homologous topological slice disk.
We finish with some questions and suggestions for
further directions in \autoref{sec:questions}.

\subsection{Conventions}

In the literature, 
properly embedded slice disks in a $4$-manifold $X$ are 
often assumed to be null-homologous
in $H_{2}(X, \partial X)$.
We will make this extra assumption on homology only
in \autoref{sec:universal_slicing} when discussing the ``universal slicings''.
For the first part \emph{deep slice} and \emph{shallow slice}
will describe the existence of
a embedded disks with the relevant properties
without conditions on the homology class.

Starting from
an $n$-manifold $M^{n}$ without boundary,
we obtain a \emph{punctured $M$}
(more precisely a \emph{bounded punctured $M$})
by removing
a small open $n$-ball $M^{\circ} \coloneqq M \setminus \interior \disk{n}$,
which yields a
manifold with boundary
$\partial M^{\circ} = \sphere{n-1}$.
Observe that a punctured $M$ is the same
as a connected sum
$M^{\circ} \cong M \# \disk{n}$
with a $n$-ball.

\subsection{Acknowledgments}

The authors would like to thank
Anthony Conway, Rob Kirby, Mark Powell, 
Arunima Ray and Peter Teichner for helpful conversations,
their encouragement and guidance.
BR would like to thank Thorben Kastenholz
for asking about the decidability of the embedded genus problem
in 4-manifolds, which motivated \autoref{sec:universal_slicing}.
We are especially grateful to
Akira Yasuhara for pointing us to related literature
and the work of Suzuki.
We are also especially grateful to the anonymous referee for numerous helpful
suggestions that improved the exposition and a suggestion that
lead to our subsequent work relating the techniques in this paper
to the embeddings of traces of surgeries on links.  
The Max Planck Institute for Mathematics in Bonn
supported us financially and with a welcoming research environment.

\section{Nonexistence of deep slice knots}
\label{sec:nonexistence_deep_slice}

\noindent For starters, we have:
\begin{proposition} \label{1-handles}
There are no deep slice knots in
$\natural^k \sphere{1} \times \ball{3}$.
\end{proposition}

\begin{proof}
Let
$K \subset \#^k \sphere{1} \times \sphere{2} 
= \partial(\natural^k \sphere{1} \times \ball{3})$
such that $K$ is slice in $\natural^k \sphere{1} \times \ball{3}$.
Then, thinking of $\natural^k \sphere{1} \times \ball{3}$
as a wedge of $k$ copies of $\sphere{1}$ thickened to be 4-dimensional,
if $D$ is any slice disk for $K$ we can isotope $D$ such that
it does not intersect a one-dimensional wedge of circles that
$\natural^k \sphere{1} \times \ball{3}$ deformation retracts onto.
Therefore, $D$ can be isotoped to be contained in
a collar neighborhood of the boundary $\#^k \sphere{1} \times \sphere{2}$
and thus $K$ is shallow slice.   
\end{proof}

The following might be a surprise,
as one could expect that additional topology
in a 3-manifold $M^3$
creates more room for concordances:
\begin{proposition}[{Special case of \cite[Prop. 2.9]{nagel2019smooth}}]
    \label{prop:local_knot_slice}
    If a local knot $K \subset \ball{3} \subset M^{3}$
    is null concordant in $M^{3} \times \interval$,
    then $K$ is null concordant in $\sphere{3} \times \interval$. 
\end{proposition}

\begin{proof}[Proof sketch]
    Let $D$ be a properly embedded disk in
    $M \times \interval$ with boundary $K$
    and let $\widetilde{M}$ be the universal cover of $M$.
    Then $D$ lifts to a properly
    embedded disk $\widetilde{D} \subset \widetilde{M} \times \interval$.
    Further, since $K$ is contained in a 3-ball $B$,
    all of the lifts of $K$ to $\widetilde{M}$ are just copies of $K$,
    and therefore, the boundary of $\widetilde{D}$
    is a copy of $K$, considered inside of $\widetilde{M}$.
    As a consequence of geometrization \cite{perelman2003finite},
    we know that every universal cover
    of a punctured compact 3-manifold smoothly embeds into $\sphere{3}$,
    as was observed in \cite[Lem.\ 2.11]{boden2017concordance}.
    It follows then that there is an embedding
    $\widetilde{M} \times \interval \hookrightarrow \sphere{3} \times \interval$.
    But then the image of $\widetilde{D}$
    under this embedding shows that $K$ bounds a disk in $\sphere{3} \times \interval$. 
\end{proof}

We have added a proof of this proposition here to
highlight that this lifting argument 
breaks down in the case of higher genus
surfaces if their inclusion
induces a nontrivial map on fundamental groups.
If $K$ bounds a genus $g$ surface 
with one boundary component $\Sigma_{g, 1}$
in $M \times \interval$, we can only
lift this to the universal cover (and subsequently find
a genus $g$ surface for $K$ in $\sphere{3}$ via this method)
under the condition that the inclusion
of $\Sigma_{g, 1}$ in $M \times \interval$ is $\pi_{1}$-trivial.
So this argument does not work if the surface
really ``uses the extra topology of $M$''.

\begin{figure}[ht]
    \begin{subfigure}{.495\textwidth}
      \centering
      \begin{overpic}[width=\textwidth]{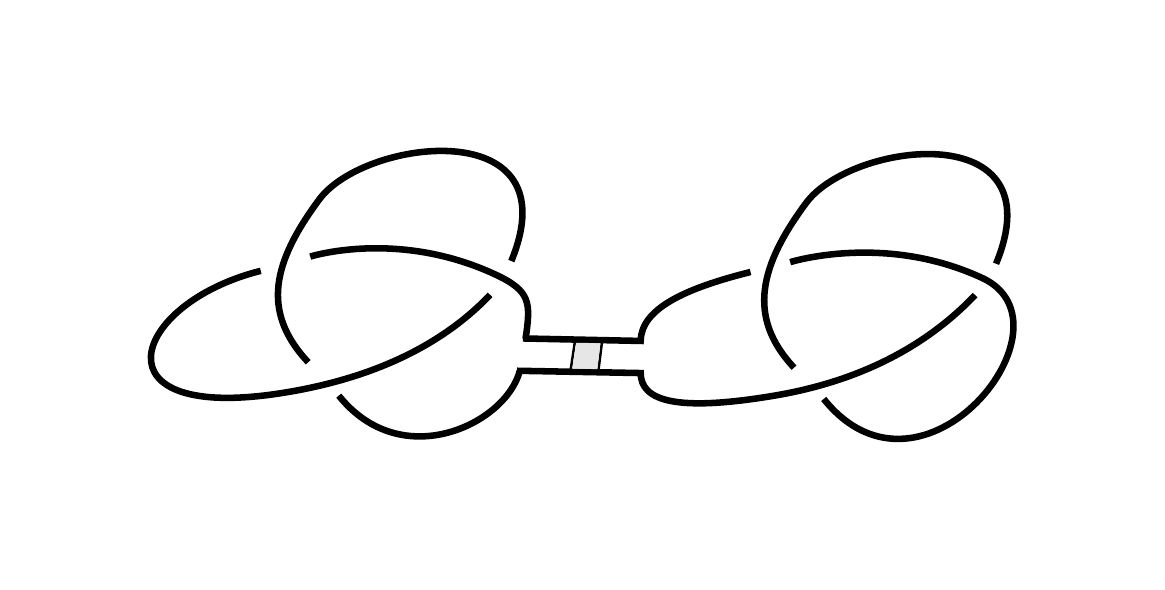}
        \label{fig:torus_movie_frame_1}
      \end{overpic}
      \caption{
            Saddle move to separate the summands
            of $K \# K$.
        }
    \end{subfigure}
    \begin{subfigure}{.495\textwidth}
      \centering
      \begin{overpic}[width=\textwidth]{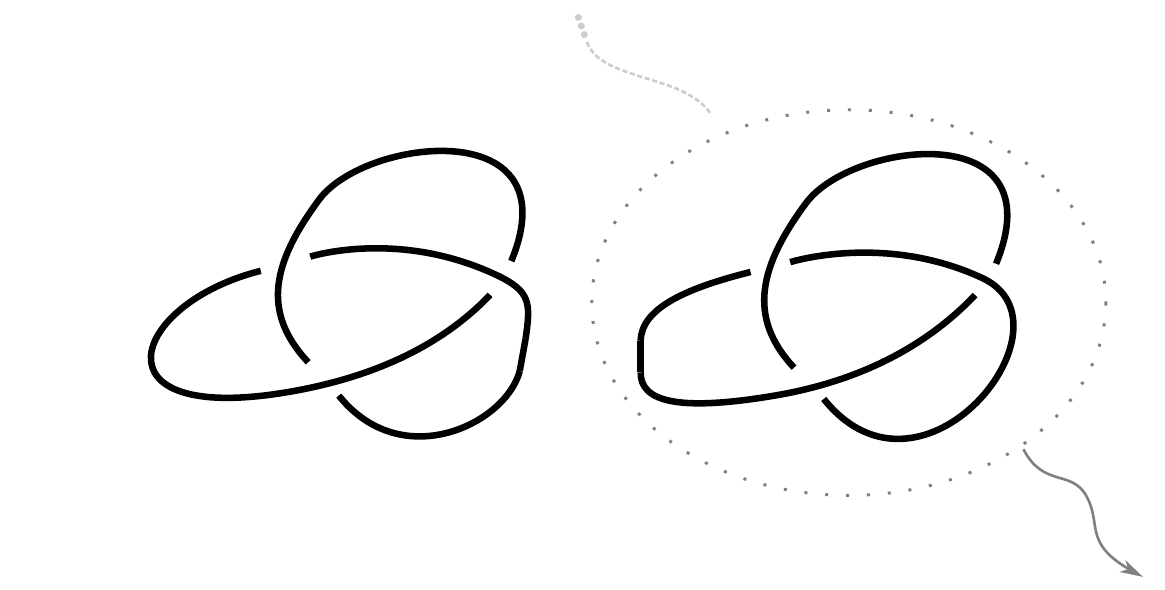}
        \label{fig:torus_movie_frame_2}
      \end{overpic}
      \caption{
            One of the summands travels around an orientation
            reversing loop in $M$.
        }
    \end{subfigure}
    \begin{subfigure}{.495\textwidth}
      \centering
      \begin{overpic}[width=\textwidth]{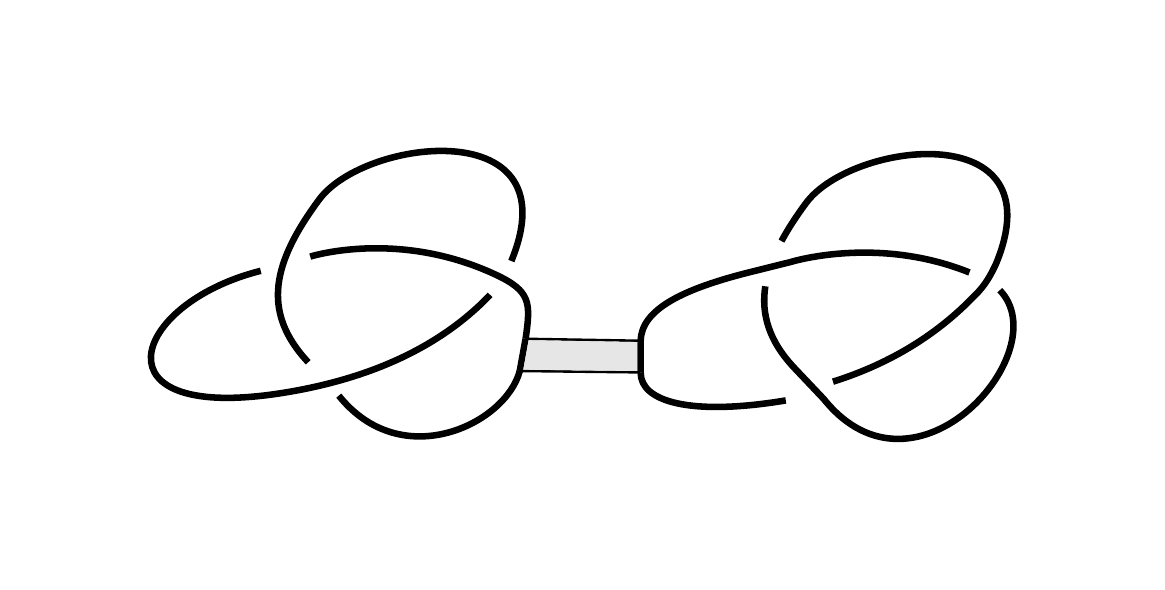}
        \label{fig:torus_movie_frame_3}
      \end{overpic}
      \caption{
            It returns mirrored, now add a fusion band.
        }
    \end{subfigure}
    \begin{subfigure}{.495\textwidth}
      \centering
      \begin{overpic}[width=\textwidth]{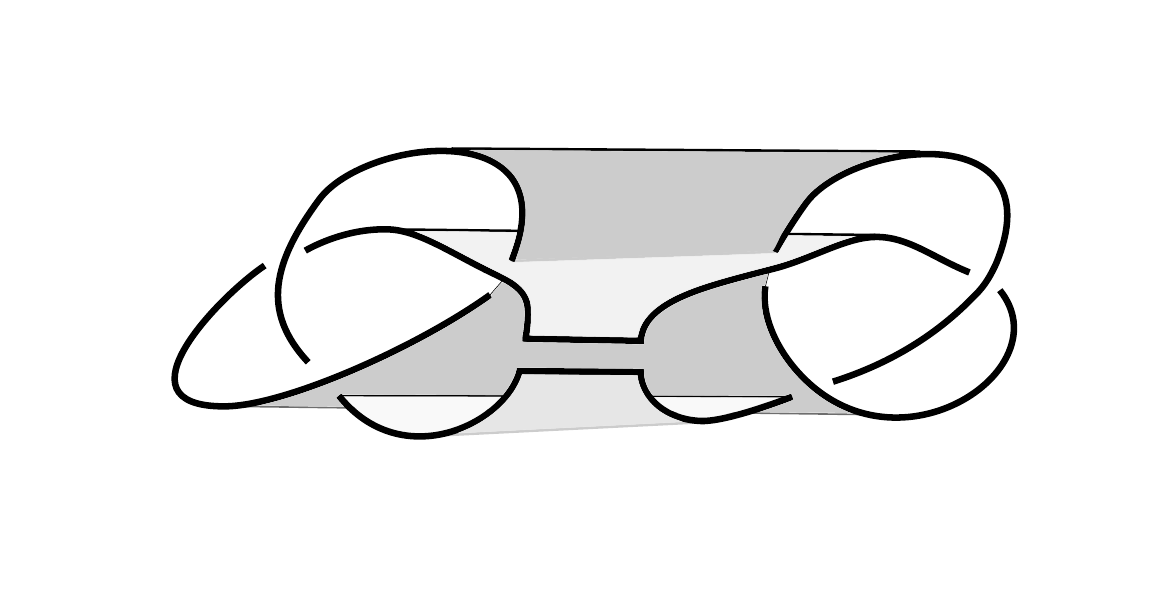}
        \label{fig:fig:torus_movie_frame_4}
      \end{overpic}
      \caption{
            Finish off the movie with the standard
            ribbon disk for $K \# -K$.
        }
    \end{subfigure}
    \caption{
        Four frames of the movie of a properly embedded
        punctured torus in $M^3 \times \interval$ with boundary
        $K \# K \subset M \times \{ 0 \}$,
        where $M$ is a non-orientable 3-manifold.
        }
    \label{fig:torus_movie}
\end{figure}

\begin{example}
    Take a non-orientable
    3-manifold $M$ containing the connected
    sum $K \# K$ of two copies of
    a local invertible knot $K$ with
    smooth 4-ball genus $g^{4}(K \# K) \ge 2$.
    As an explicit example,
    $K$ a left-handed trefoil will work,
    and we illustrate the following
    in \autoref{fig:torus_movie}.
    Describe an embedded torus in
    $M \times \interval$ with the motion picture method:
    Use the connected sum band to split the sum with a saddle.
    Then let one of the summands travel around an orientation
    reversing loop in $M$
    while leaving the other one fixed.
    The summand traveling around the loop
    was reflected in the process and since it is invertible
    it is isotopic to
    $-K = r \overline{K}$
    in a 3-ball neighborhood in $M$.
    Fusing the summands back together
    along a connected sum band
    we now obtain
    $K \# -K$ as a local knot.
    Finally cap this off with the usual ribbon disk
    for the connected sum of a knot with
    its concordance inverse.  Therefore $g^{M \times I}(K \# K) \leq 1$ and by
    \autoref{prop:local_knot_slice} in fact $g^{M \times I}(K \# K) = 1$.
    The 4-ball genus $g^{4}(K \# K) \ge 2$
    of this example is strictly
    larger than its
    \emph{$(M^{3} \times \interval)$-4-genus},
    which we define as
    \[
        g^{M \times \interval}(J)
        \coloneqq
        \min 
        \{ 
            g \mid
            \exists \text{ smooth proper embedding }
            \Sigma_{g, 1} \hookrightarrow M \times I
            \text{ with }
            \partial \Sigma_{g, 1} = J \subset M \times \{ 0 \}
        \}
    \]
    Observe that in this notation
    the usual 4-ball genus is
    $g^{4} = g^{\sphere{3} \times \interval}$
    and we can rephrase \autoref{prop:local_knot_slice} as
    $g^{M \times \interval}(K) = 0$ implies
    $g^4(K) = 0$
    for local knots $K$.
    Similar notions of 4-genera
    were introduced in Celoria's investigation
    of almost-concordance \cite[Def.\ 12]{celoria2018concordance}.
\end{example}

It would be interesting to find an example 
of an orientable $3$-manifold $M^3$ where the
$g^{M \times \interval}(K)$ genus of
some local knot $K \subset \disk{3} \subset M$
is strictly smaller than
the $4$-ball genus $g^{4}(K)$, or prove that no such $M$ exists.
Local $K$ satisfy $g^{M \times \interval}(K) \le g^4(K)$
as cobordisms in $\sphere{3} \times \interval$
can be embedded into $M \times \interval$.
Because of \autoref{prop:local_knot_slice}
an example where these values differ can only appear
for $g^{4}(K) \ge 2$.
Moreover, as we will see in \autoref{thm:embedding_4-sphere_shallow}
such an $M$ would necessarily not embed in $S^4$.  
Another special case is treated in
\cite[Thm.\ 2.5]{Davis_2018} where 
a handle cancellation argument shows that
there is no difference
for local knots in $M = \sphere{1} \times \sphere{2}$,
that is the equality 
$g^{\sphere{1} \times \sphere{2} \times \interval} = g^{4}$ holds
(and also analogous statements for $\#^{k} \sphere{1} \times \sphere{2}$).
Topological concordance in $\sphere{1} \times \sphere{2} \times \interval$
is investigated in \cite{friedl2019satellites}.

We now give a criterion that shows that certain 4-manifolds
have no local deep slice knots in the boundary.
This idea is also contained in \cite[Thm.\ 0]{suzuki1969localknots}
and its variants.

\begin{proposition} 
    \label{thm:embedding_4-sphere_shallow}
    Let $X^4$ be a compact 4-manifold with a 
    local knot $\gamma \subset \ball{3} \subset \partial X$
    that is slice in $X$.
    If there is a cover of $X$ which can be smoothly embedded into $\sphere{4}$,
    then $\gamma \subset \ball{3} \hookrightarrow \sphere{3} = \partial \ball{4}$
    is slice in $\ball{4}$.
    Hence, $\gamma$ is shallow slice in $X$.    
\end{proposition}

\begin{proof}
    Let $\widetilde{X}$ be a cover of $X$
    with an embedding $\widetilde{X} \subset \sphere{4}$ into $\sphere{4}$
    and let $\widetilde{D}$ be a lift of a slice disk for $\gamma$ to
    $\widetilde{X}$
    with $\widetilde{\gamma} = \partial \widetilde{D}$.
    Note that the knot $\widetilde{\gamma}$ is the same as $\gamma$,
    since $\gamma$ is contained in a 3-ball and the only covers of a 3-ball
    are disjoint unions of 3-balls. Puncture $S^4$ by removing a small ball $B$ close to $\widetilde{\gamma}$ and such that $\widetilde{\gamma}$ can be connected by an annulus  disjoint from $\tilde{X}$ to $\partial B$ and such that the other end of the annulus is (the mirror image of) $K \subset \partial B$.
    Then since $S^4 - \interior{B} \cong B^4$, 
    the annulus together with $\widetilde{D}$ show that $K$ is slice in
    the $B^4$ which is the complement of the small ball.    
    Therefore $\gamma$ is shallow slice in $X$.  
\end{proof}

As an example, \autoref{thm:embedding_4-sphere_shallow}
implies that $\natural^k \sphere{2} \times \disk{2}$
contains no deep slice local knots,
since these manifolds can all be embedded in $\sphere{4}$.
However, these manifolds all contain deep slice knots, necessarily non-local, 
as will be seen shortly.
Additionally, we have:
\begin{corollary}
    \label{cor:universal_cover_S4_no_deep_slice}
    Suppose that $X$ is a closed smooth 4-manifold with universal cover $\mathbb{R}^4$ or $\sphere{4}$,
    and let $X^{\circ}$ denote the punctured version.
    Then $X^{\circ}$ has no deep slice knots.  
\end{corollary}

\section{Existence of deep slice knots}
\label{sec:existence_deep_slice}

A \emph{2-handlebody} is a 4-manifold whose handle decomposition
contains one 0-handle, some nonzero number of 2-handles
and no handles of any other index.
Examples of this are knot traces,
where a single $2$-handle is attached along a framed knot
to the 4-ball.
In this section, we prove:
\begin{theorem} 
    \label{existence}
    Every 2-handlebody $X$ contains a null-homotopic
    deep slice knot in its boundary.  
\end{theorem}

\begin{remark}
    For the special case of the $2$-handlebody
    $\disk{2} \times \sphere{2}$
    the existence of such knots was
    already observed in
    \cite[Thm.\ B]{Davis_2018}
    (here only winding number $w = 0$ gives null homotopic knots).
    Furthermore the authors construct an infinite family of slice knots
    which are pairwise different in topological concordance in
    a collar of the boundary.
\end{remark}

\autoref{existence} breaks up naturally into two cases depending on whether
the boundary has nontrivial $\pi_1$ or not (i.e. if it is or is not $S^3$).
In the case where $\pi_1(\partial X) \neq 1$,
there is a concordance invariant for knots in arbitrary 3-manifolds,
closely related to the Wall self-intersection number
(see \cite{wall1999surgery}, \cite{freedman1990topology},
and \cite{schneiderman2003algebraic}),
that will allow us to show that some obviously slice knots are not shallow slice.
In the case where $\pi_1 (\partial X)$ is trivial,
and therefore by the 
3-dimensional Poincar\'e conjecture \cite{perelman2003finite}
$\partial X = \sphere{3}$, the Wall self-intersection number is of no use.
However, in this case, the consideration of whether a knot
that is slice in $X$ is deep slice in $X$ is related
to the existence of spheres representing various homology
classes in the manifold obtained by closing $X$ off with a 4-handle.   

\begin{remark}
    If there was a direct proof
    that every closed homotopy 3-sphere
    smoothly bounds a contractible 4-manifold,
    then we would not need to
    invoke the $3$-dimensional Poincar\'e conjecture.
\end{remark}

Following \cite{yildiz2018note} and \cite{schneiderman2003algebraic},
we briefly introduce the Wall self-intersection number
in the setting that we will be working in, and state some of its basic properties.
Let $Y^3$ be a closed oriented 3-manifold and let
$\gamma \colon \sphere{1} \hookrightarrow Y$ be a knot in $Y$.
Let $\mathcal{C}_\gamma(Y)$ denote the set of concordance classes
of oriented knots in $Y$ that are freely-homotopic to $\gamma$.
In particular $\mathcal{C}_{U}(Y)$
denotes the set of concordance classes of
oriented null-homotopic knots in $Y$,
where we write $U$ for the local unknot in $Y$.
Given an  oriented null-homotopic knot $K \subset Y$,
by transversality there exists an oriented immersed disk
$D$ in $Y \times \interval$ with boundary
$K \subset Y \times \{ 0 \} = Y$
that has only double points of self-intersection.
Let $\star \in Y$ denote a basepoint which we implicitly use for
$\pi_1(Y) = \pi_1(Y \times \interval)$ throughout.
Choose an arc, which we will call a whisker, from $\star$ to $D$.
For each double point of self-intersection $p \in D$
choose a numbering of the two sheets of $D$ that intersect at $p$.
Then let $g_p \in \pi_1(Y)$ be the homotopy class of the loop in $Y \times \interval$
obtained by starting at $\star$, taking the whisker to $D$,
taking a path to $p$ going in on the first sheet,
taking a path back to where the whisker meets $D$ that leaves $p$ on the second sheet,
and then returning to $\star$ using the whisker.
Note that changing the order of the two sheets would transform $g_p$ to $g_p^{-1}$.
Also, since $K$ and $Y$ are oriented,
$D$ and $Y \times \interval$ obtain orientations
with the convention that $K \subset Y \times \{ 0 \} = Y$,
and therefore, for every self-intersection point $p \in D$,
there is an associated sign which we will denote by $\operatorname{sign}(p)$.  

\noindent Let
\[
    \widetilde{\Lambda}
    \coloneqq
    \frac{   \mathbb{Z}[\pi_1(Y)]  }
	{  \langle \{g - g^{-1} \mid g \in \pi_1(Y)\} \rangle \oplus \mathbb{Z}[1]  }
\]
were the quotient is a quotient as abelian groups.
The \emph{Wall self-intersection number} of $K$ is defined to be 
\[
    \mu(K) = \sum_p \operatorname{sign}(p) \cdot g_p \in \widetilde{\Lambda}
\]
See \cite{schneiderman2003algebraic}
for a proof that it is independent of the choice of $D$, the choice of whisker,
and the choice of orderings of the sheets of $D$ around the double points.
Further, $\mu$ is a concordance invariant in $Y \times \interval$,
and therefore defines a map:
\[
    \mu \colon \mathcal{C}_{U}(Y) \to \widetilde{\Lambda}
\]
Notice that if $g \in \pi_1(Y)$ and $g \neq 1$ then $g$ is also nonzero in $\widetilde{\Lambda}$.  

\begin{figure}
    \centering
    \begin{overpic}[width=0.9\textwidth]{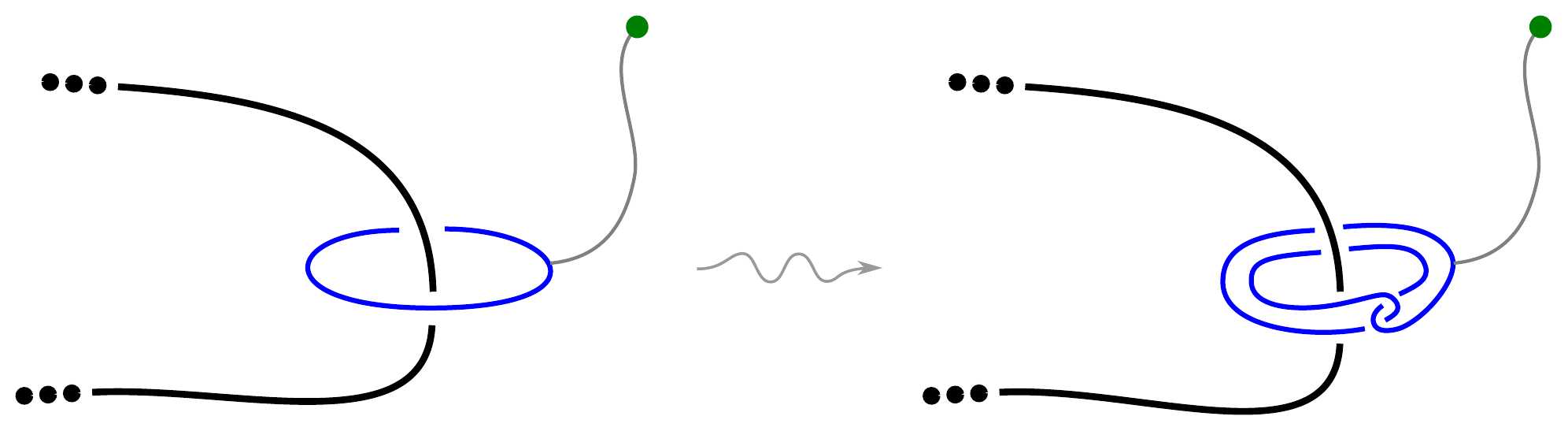}
        \put(5, 18){$L_{i}$}
        \put(15, 10){\color{blue} $\gamma$}
        \put(36.5, 25){$\star$}
        \put(63, 10){\color{blue} $\operatorname{Wh}(\gamma) = K$}
        \put(63, 18){$L_{i}$}
        \put(94, 25){$\star$}
    \end{overpic}
    \caption{The Whitehead double of a nontrivial meridian
    $\gamma$ to one of the surgery link components is deeply slice in $X$.}
    \label{fig:meridian_Whitehead_double}
\end{figure}

\begin{figure}
    \centering
    \begin{overpic}[width=0.9\textwidth]{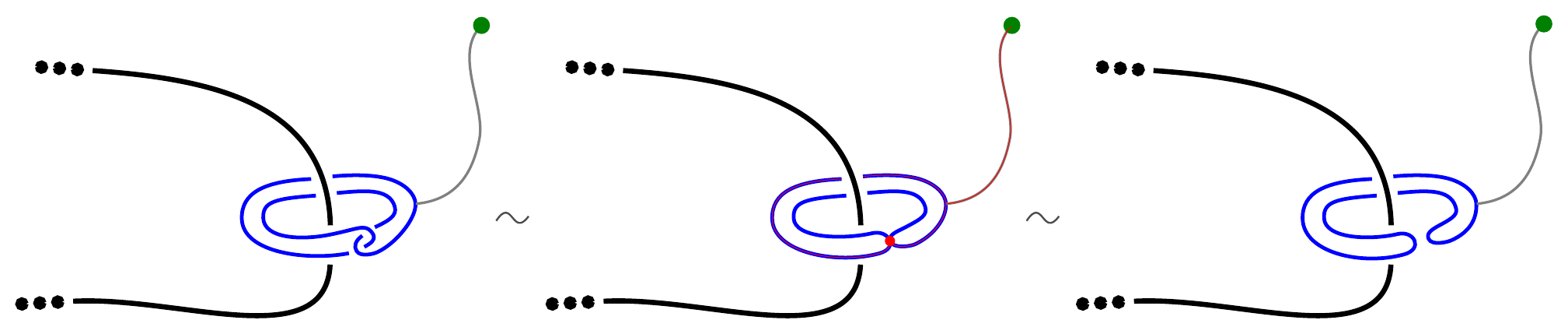}
    \end{overpic}
    \caption{
        Track of a homotopy from the Whitehead double of $\gamma$
        to the unknot giving an immersed disk with a
        single double point (red in the middle frame).
        The red double point loop based at the green basepoint
        calculates that $\mu(K) = \gamma$.
    }
    \label{fig:Whitehead_double_nullhomotopy}
\end{figure}

\begin{proof}[{Proof of \autoref{existence}, Case 1}]
We are now in position to handle \autoref{existence}
in the case where $\pi_1(\partial X) \neq 1$.
Now $X$ is described by attaching 2-handles to $\disk{4}$
along some framed link $L \subset \partial \disk{4}$.
Since $\pi_1(\partial X)$ is 
(normally) generated by the meridians of $L$
and $\pi_1(\partial X) \neq 1$, there is some meridian $\gamma$ of
$L$ that is nontrivial in $\pi_1(\partial X)$.
Notice that if we are given a 2-handlebody described by a framed link $L$
and $K$ is a knot in the boundary of the 2-handlebody that is shown in the
framed link diagram as an unknot (possibly linked with $L$),
then $K$ is slice in the 2-handlebody -- just forget all of the other
2-handles and take an unknotting disk
whose interior has been pushed into the 0-handle. 
Now, take $K$ in $\partial X$ to be a
Whitehead double of $\gamma$ as in \autoref{fig:meridian_Whitehead_double},
which is a null-homotopic knot in the boundary. 
By the previous observation, since $K$ is unknotted in the boundary of the 0-handle,
$K$ is slice in $X$.
Additionally, one computes that $\mu(K) = \gamma \neq 1 \in \widetilde{\Lambda}$,
for example using the null-homotopy in \autoref{fig:Whitehead_double_nullhomotopy}.
Therefore, $K$ is not null-concordant in $\partial X$, so $K$ is deep slice in $X$.
\end{proof}

Notice that if $\pi_1(\partial X) = 1$,
then $\mu$ is of no use since $\widetilde{\Lambda} = 0$.
Now assume that $\pi_1(\partial X) = 1$ so that $\partial X = \sphere{3}$.
Again $X$ is obtained by attaching 2-handles to some framed link $L$.
Let $\widehat{X}$ denote the closed 4-manifold
obtained by closing off $X$ with a 4-handle.
We will need a lemma on surfaces in 2-handlebodies,
whose statement is
standard and could alternatively be concluded from
the KSS-normal form for surfaces
as in \cite[{Thm.\ 3.2.7}]{kamada2017surface} and \cite{MR672939}.
\begin{lemma}
    \label{lem:disk_in_the_bottom}
    Let $X$ be a closed smooth 4-manifold with a handle decomposition consisting
    of only 0-, 2-, and 4-handles,
    with exactly one 0-handle and one 4-handle.
    Every element of $H_2(X;\mathbb{Z})$ can be represented
    by a smooth closed orientable surface whose
    intersection with the union of the 0- and 2-handles of $X$
    is a single disk.
\end{lemma}

\begin{proof}	
Let $X_{\le 2}$ denote the union of the 0- and 2-handles of $X$,
so that $X = X_{\le 2} \cup \ball{4}$.
For every 2-handle $h_i$, there is an element $H_2(X; \mathbb{Z})$
obtained by taking the co-core disk $D_i$ for $h_i$
and capping it off
with an orientable surface in the 4-handle. 
Let $\{F_i\}$ denote a choice of these surfaces, one for each 2-handle.
These surfaces form a basis for $H_2(X;\mathbb{Z})$ and note
that each has the desired property that $F_i \cap X_{\le 2} = D_i$ is a disk.

Given an arbitrary element $x \in H_2(X; \mathbb{Z})$, we have
$x = a_1 [F_1] + \cdots + a_n[F_n]$ for some $a_i \in \mathbb{Z}$.
Therefore, by taking parallel copies of the $F_i$ for each summand,
we can find an immersed (possibly disconnected)
orientable surface $F'$ representing $x$,
with $F' \cap X_{\le 2}$ a union of $\sum \abs{a_{i}}$ disjoint disks.
By taking arcs in $\partial X_{\le 2}$ that connect
the different boundaries of the disks all together,
and attaching tubes to $F'$ along these arcs,
we obtain a connected orientable immersed surface $F''$ representing $x$
whose intersection with $X_{\le 2}$ is now a disk.
In particular, the tubing is done so that half of the tube is contained in $X_{\le 2}$ and the other half is in the 4-handle, and therefore $F'' \cap X_{\le 2}$ is the result of boundary summing together the disks in $F' \cap X_{\leq 2}$.

To make $F''$ into an embedded surface,
we can resolve the double points in the 4-handle,
by increasing the genus,
and arrive at a surface
representing $x$ with the desired property.  
\end{proof}

The main ingredient for the proof of the second case of
\autoref{existence} is the following theorem of Rohlin,
and in particular the corollary that follows.
Rohlin's theorem has been used in a similar way
to study slice knots in punctured connected sums
of projective spaces, for example in \cite{yasuhara1991torusknot} 
and \cite{yasuhara1992complexplane}.
\begin{theorem}[{Rohlin, \cite{rokhlin1971two}}]
	\label{Rohlin}
	Let $X$ be an oriented closed smooth 4-manifold with $H_1(X; \mathbb{Z}) = 0$.
	Let $\psi \in H_2(X; \mathbb{Z})$ be an element that is divisible by 2,
	and let $F$ be a closed oriented surface of genus $g$
	smoothly embedded in $X$ that represents $\psi$. Then
	$$
		4g \geq | \psi \cdot \psi - 2\sigma(X) | - 2b_2(X)
	$$
\end{theorem}

\begin{corollary}
    \label{cor:no_sphere}
    Let $X$ be a closed smooth 4-manifold with $H_1(X; \mathbb{Z}) = 0$, and $H_2(X; \mathbb{Z}) \neq 0$.  Then there exists a homology class $\psi \in H_2(X; \mathbb{Z})$ that cannot be represented by a smoothly embedded sphere.  
\end{corollary}

\begin{proof}[Proof of \autoref{cor:no_sphere}]
To apply \autoref{Rohlin}, we must find a homology class $\psi$
that is divisible by 2 where the right hand side
$| \psi \cdot \psi - 2\sigma(X) | - 2b_2(X) > 0$.
Since the intersection form on $X$ is unimodular,
there exists some element $\alpha$ with $\alpha \cdot \alpha \neq 0$.
From Poincar{\'e} duality together with the universal coefficient
theorem and our hypothesis that $H_{1}(X; \ZZ) =0$, we know that $H_{2}(X; \ZZ)$ 
is torsion free.
Then by taking $k$ to be a sufficiently large integer,
we can make $|(2k \alpha) \cdot (2k \alpha) - 2 \sigma(X) |$ arbitrarily large.
By taking $\psi = 2k \alpha$, the result follows.  
\end{proof}

\begin{proof}[{Proof of \autoref{existence}, Case 2}]
By \autoref{cor:no_sphere}, let $\psi \in H_2(\widehat{X}; \mathbb{Z})$
be a homology class that can not be represented by an embedded sphere.
Using \autoref{lem:disk_in_the_bottom}, let $F$ be a smooth closed orientable surface
representing $\psi$ whose intersection with $X = \widehat{X}_{\le 2}$ is a disk $D$,
as illustrated schematically in \autoref{fig:surface_schematic}.
\begin{figure}
    \centering
    \begin{overpic}[width=0.45\textwidth]{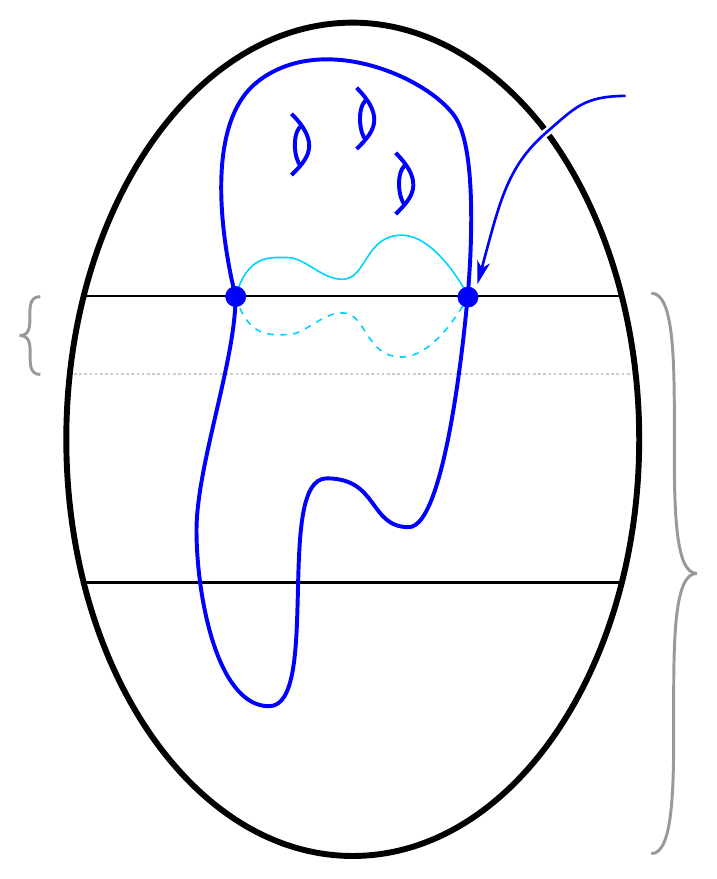}
        \definecolor{blue-green}{rgb}{0.0, 0.87, 0.87}
        \put(73, 88){\color{blue} $\partial D$}
        \put(-11.5, 60.5){\color{gray} $\sphere{3} \times \interval$}
        \put(60, 69){4-h.}
        \put(81, 33){$X = \widehat{X}_{\le 2}$}
        \put(60, 36){2-h.}
        \put(60, 28){0-h.}
        \put(36, 20){\color{blue} $D$}
    \end{overpic}
    \caption{
        Schematic of the blue surface $F$
        in the 2-handlebody, intersecting $X=$ the
        union of the 0- and 2-handles
        in a disk $D$. If $\partial D$
        was shallow slice (dashed light blue) in $X$,
        disk $D$ union the shallow slice disk
        flipped into the 4-handle (solid light blue)
        would be an impossible sphere representative
        of the homology class of $F$.
    }
    \label{fig:surface_schematic}
\end{figure}
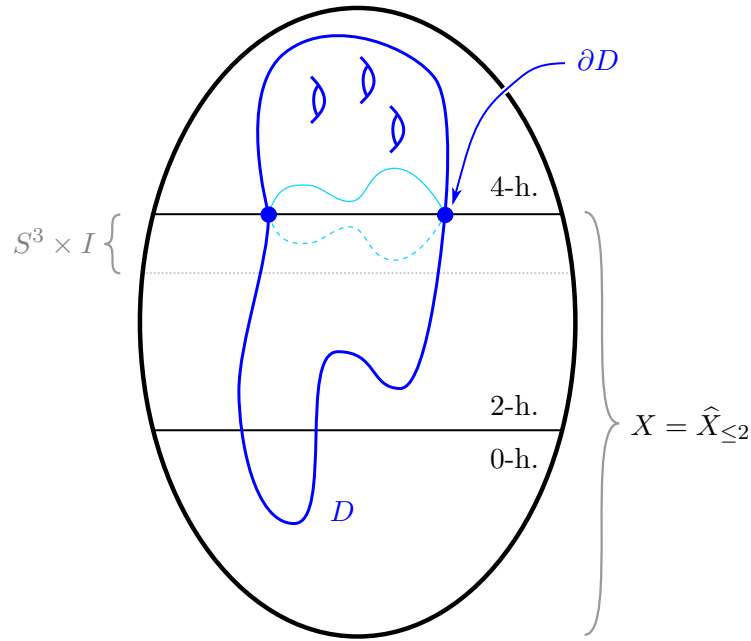
Then $\partial D \subset \partial X$ is deep slice in $X$,
since otherwise the surface obtained by intersecting $F$
with the 4-handle could be replaced with a disk
without altering the homology class,
violating the assumption that $\psi$ cannot be represented
by an embedded sphere.
To see that the homology class is not altered,
observe that in any 2-handlebody the homology class of a surface
is determined by how it intersects the 0- and 2-handles. 
Also observe that in this case the deep slice knot
$\partial D \subset \partial X = \sphere{3}$
is local.
This concludes the proof of \autoref{existence}.
\end{proof}

\section{Universal slicing manifolds do not exist}
\label{sec:universal_slicing}

The Norman-Suzuki trick
\cite[Cor.\ 3]{MR246309}, \cite[Thm.\ 1]{suzuki1969localknots}
can be used to show
that any knot $K \subset \sphere{3}$ bounds
a properly embedded disk in a punctured
$\sphere{2} \times \sphere{2}$:
The track of a null-homotopy of $K$ in $\disk{4}$
can be placed in the punctured $\sphere{2} \times \sphere{2}$
which gives a disk
that we can assume to be a generic immersion,
missing $\sphere{2} \vee \sphere{2} \subset (\sphere{2} \times \sphere{2})^{\circ}$,
and with a finite number
of double points.
By tubing into the spheres
$\sphere{2} \times \{ \textrm{pt} \}, \{ \textrm{pt} \} \times \sphere{2}$
we can remove all the intersections
-- but observe that
this changes the homology class of the disk.

\begin{proposition}
    \label{prop:all_slice}
	Let $M^3$ be a closed orientable 3-manifold.
	There exists a compact orientable 4-manifold $X^4$ constructed
	with only a 0-handle and 2-handles, with $\partial X = M$
	such that every knot in $M$ is slice in $X$.
\end{proposition}

\begin{proof}
	Start by taking any compact 4-manifold $X'$
	with only 0-, 2-handles
	and boundary $M$
	and let $X = X' \# (\sphere{2} \times \sphere{2})$.
	Let $K \subset M = \partial X$ be a knot.
	Since $X'$ and $X$ are simply connected,
	$K$ bounds an immersed disk which we can assume
	lives completely in the $X'$-summand of the connected sum.
	Now the Norman-Suzuki trick works to remove intersection points
	of the immersion
	by tubing into the coordinate spheres
	of the $\sphere{2} \times \sphere{2}$-summand.
\end{proof}

\begin{remark}
    In contrast to the homologically nontrivial
    disks constructed in the Norman-Suzuki trick,
    a knot is slice via a null-homologous disk in
    some connected sum $\#^{n} \sphere{2} \times \sphere{2}$
    if and only if its Arf-invariant is zero.
    $\operatorname{Arf} K = 0$ implies that the knot is band-pass
    equivalent to the unknot, and a band pass
    can be realized by sliding the (oppositely oriented) strands
    of a pair of bands over the coordinate spheres
    in a $\sphere{2} \times \sphere{2}$ factor.
    Conway-Nagel \cite{conway2020stably}
    defined and studied the minimal number of summands
    needed 
    to find a disk 
    in a punctured $\#^{n} \sphere{2} \times \sphere{2}$.
\end{remark}

\noindent \textbf{Convention:} From now until the
end of this section, properly embedded
slice disks $\Delta^{2} \subset X^{4}$ in a $4$-manifold
are always required
to be \textbf{null-homologous}.
We will still add the qualifier ``null-homologous''
in the statements to emphasize this.
Since our obstructions work in the topologically locally
flat category, we will formulate everything in this more general setting.

\begin{definition}
    A knot $K \subset \sphere{3}$ is
    \emph{(topologically/smoothly)
    null-homologous slice in the (topological/smooth)
    $4$-manifold $X^{4}$} with
    $\partial X = \sphere{3}$,
    if $K = \partial \Delta$, where
    $\Delta^{2} \subset X$ is a (locally flat/smooth) properly embedded
    disk such that
    $[\Delta, \partial \Delta] = 0 \in H_{2}(X, \partial X)$.
\end{definition}

One way of studying if a knot $K$ is slice in $\disk{4}$
is to approximate $D^4$ by varying the $4$-manifold $X$.
By restricting the intersection form
and looking at simply-connected 4-manifolds $X$
this gives rise to various filtrations of the
knot concordance group
(notably the $(n)$-solvable filtration $\mathcal{F}_{n}$
of Cochran-Orr-Teichner \cite{COT2003}
and the positive and negative variants
$\mathcal{P}_{n}, \mathcal{N}_{n}$ \cite{MR3109864}).

We say that the properly embedded disk
$\Delta$ is \emph{null-homologous}
if its fundamental class
$[\Delta, \partial \Delta] \in H_{2}(X, \partial X)$
is zero.
Since by Poincar{\'e} duality the intersection pairing 
$H_{2}(X) \otimes_{\ZZ} H_{2}(X, \partial X) \xrightarrow{\pitchfork} \ZZ$
is non-degenerate, 
a null-homologous disk is characterized
by the property that it intersects all closed second homology classes
algebraically zero times.
For slicing in arbitrary $4$-manifolds,
we here restrict to null-homologous disks
to exclude constructions
as in the Norman-Suzuki trick.

For every fixed knot $K \subset \sphere{3}$,
there is a $4$-manifold in which $K$ is null-homologically slice.
Norman \cite[Thm.\ 4]{MR246309} already
observes that it is possible
to take as the 4-manifold a punctured connected sum
of the twisted 2-sphere bundles $\sphere{2} \widetilde{\times} \sphere{2}$.
Similarly, \cite[Lem.\ 3.4]{cochran1986unknotting} discuss that
for any knot $K \subset \sphere{3}$
there are numbers $p, q \in \mathbb{N}$
such that $K$ is null-homologous slice
in the punctured connected sum 
$\#^{p} \CP{2} \#^{q} \overline{\CP{2}}$ of
complex projective planes.
The argument starts with a sequence of positive and negative
crossing changes leading from $K$ to the unknot, 
and then realizes say a positive crossing change
by sliding a pair of oppositely
oriented strands over the $\CP{1}$
in a projective plane summand. The track of this
isotopy, together with a disk bounding the final unknot
gives a motion picture of a null-homologous slice disk.
Since both positive and negative crossing changes might be
necessary, it is important that both orientations
$\CP{2}, \overline{\CP{2}}$ are allowed to appear in the connected sum.

In view of $(\sphere{2} \times \sphere{2})^{\circ}$
where every knot in the boundary bounds a disk
(which is rarely null-homologous)
and $(\#^{p} \CP{2} \#^{q} \overline{\CP{2}})^{\circ}$,
in which we find plenty of null-homologous disks
(but only know how many summands $p, q$ we need
after fixing a knot on the boundary)
a natural question concerns the
existence of a
\emph{universal slicing}
manifold.
Is there a fixed compact, smooth, oriented
$4$-manifold $V^{4}$ with $\partial V = \sphere{3}$
such that any knot $K \subset \sphere{3}$
is slice in $V$ via a null-homologous disk?
It turns out that a signature estimate shows such
a universal solution cannot exist.
\begin{theorem}
	\label{thm:no_universal_slicing}
	Any compact oriented $4$-manifold $V^{4}$ 
	with $\partial V = \sphere{3}$
	contains a knot
	in its boundary
	that is not 
	topologically null-homologous slice in $V$.
\end{theorem}

\begin{remark}
    If we drop the assumption that $V$ should be compact,
    a punctured infinite connected sum of projective planes
    does the job:
    \begin{equation*}
       \disk{4} \#^{\infty} (\CP{2} \# \overline{\CP{2}})
    \end{equation*}
    For any fixed knot on the
    boundary there
    is a compact slice disk in a finite stage
    \begin{equation*}
        \disk{4} \#^{k} \CP{2} \#^{l} \overline{\CP{2}} \# \disk{4}
        \subset
        \disk{4} \#^{\infty} (\CP{2} \# \overline{\CP{2}}).
    \end{equation*}
\end{remark}

The remainder of this section is concerned with
a proof of \autoref{thm:no_universal_slicing}.
As preparation, let us specialize a result
\cite[Thm.\ 3.8]{conway2020stably},
which is a generalization of the Murasugi-Tristram inequality
for links bounding surfaces in $4$-manifolds,
to the case of knots.
Here $\sigma_{\omega}(K)$ is the
\emph{Levine-Tristram signature} of the knot $K$,
defined as the signature
of the hermitian matrix
$(1-\omega)V + (1 - \overline{\omega})V^T$,
where $V$ is a Seifert matrix of $K$
and $\omega$ a unit complex number not equal to $1$.
References for this signature include
\cite{MR253348}, \cite{tristram1969some}
and the recent survey \cite{conway2019levine}.
The following inequality only holds for specific values of $\omega$,
and will adopt the notation
$\sphere{1}_{!}$ for unit complex numbers
$\omega \in \sphere{1} - \{ 1 \}$
which do not appear as a zero of an integral Laurent polynomial
$p \in \laurent$ with $p(1) = \pm 1$.
\begin{theorem}[{\cite[Special case of Thm.\ 3.8]{conway2020stably}}]
    Let $X$ be a closed oriented topological $4$-manifold
    with $H_{1}(X; \ZZ) = 0$.
    If $\Sigma \subset (\sphere{3} \times \interval) \# X$
    is a null-homologous (topological)
    cobordism between two knots
    $K \subset \sphere{3} \times \{ 0 \}$ and
    $-K' \subset - (\sphere{3} \times \{ 1 \})$,
    each contained
    in one of the two boundary component $\sphere{3}$'s
    of $(\sphere{3} \times \interval) \# X$, then
    \begin{equation*}
         \abs{ \sigma_{K'}(\omega) - \sigma_{K}(\omega) + \sign(X) }
        - \chi(X) + 2
        \le
        - \chi(\Sigma)
    \end{equation*}
    for all $\omega \in \sphere{1}_{!}$.
\end{theorem}


\noindent For $K \subset \partial X^{\circ}$ which is
null-homologous slice in $X$ and $\Sigma$ an annulus,
we can further simplify:
\begin{corollary}
    \label{cor:MT_general_manifold}
    Let $X$ be a closed topological $4$-manifold with
    $H_{1}(X; \ZZ) = 0$.
    If the knot $K \subset \sphere{3}$ is topologically
    null-homologous slice in $X^{\circ}$ then
    for $\omega \in \sphere{1}_{!}$ we have
    \begin{equation*}
        \abs{ \sigma_{K}(\omega) + \sign(X) }
        - \chi(X) + 2
        \le
        0
    \end{equation*}
\end{corollary}

To prove \autoref{thm:no_universal_slicing}
it will be enough to obstruct the sliceness of a single knot
in the boundary.
The strategy is to use surgery to trivialize $H_{1}$,
then pick the knot $K$ in the original manifold boundary
and arrive at a contradiction to \autoref{cor:MT_general_manifold}
in the surgered manifold if $K$ was null-homologous slice.

\begin{proof}[Proof of \autoref{thm:no_universal_slicing}]
	Let $V$ be a compact topological $4$-manifold
	with boundary $\sphere{3}$, we want to
	find a knot in its boundary which is not slice.
	Pick a set of disjointly embedded loops
	$\gamma_{1}, \ldots, \gamma_{l}$ in $V$ whose homology classes
	generate $H_{1}(V)$.
	If $V$ already satisfies $H_{1}(V) = 0$,
	set $l=0$ for the remainder of the proof and
	omit the surgery altogether.
	Let $K$ be a knot in $\sphere{3}$ whose signature
	(at the unit complex number $\omega = -1$)
	satisfies
	\begin{equation*}
		\abs{\sigma_{K}(-1)} \ge \abs{\sign(V)} + \abs{\chi(V)} + 2l.
	\end{equation*}
	Note that
	the constant on the right hand side only depends
	on the signature, Euler characteristic,
	and number of generators of $H_{1}(V)$,
	and not on the knot $K$.
    For example, since signature is additive
    under connected sum,
    the self-sum $K_{n} = \#^{n} K$ with $n$ large enough
    has arbitrarily high signature at $\omega = -1$
    if we start with a $K$
    that has positive signature $\sigma_{K}(-1)$
    (for example, taking $K$ to be the left-handed trefoil knot).
	
	Suppose that $K$ is slice in $V$ via a null-homologous disk $\Delta$.
	Being null-homologous in
	the relative second homology group
	means geometrically that there is a locally flat embedded 3-manifold
	$M^3$ with boundary the slice disk $\Delta$ union
	a Seifert surface for $K$ in the boundary $\sphere{3}$,
	see \cite[Lem.\ 8.14]{lickorish1997introduction}.
	We can remove the closed components from $M$,
	what remains is a 3-manifold with nonempty boundary in $V$.
	Generically the embedded circles $\gamma_1, \ldots, \gamma_l$
	will intersect the
	$3$-manifold $M$ in points, but we can push these intersection
	points off the boundary of $M$
	via an isotopy of the curves in $V$.
	We will still keep the notation
	$\gamma_1, \ldots, \gamma_l$ for the isotoped curves
	which are now disjoint from $M$.
	Essentially, this finger move supported in a
	neighborhood of $M$ is guided by pairwise disjoint arcs in $M$
	connecting the intersections points to the boundary.
	
	Perform surgery on the loops
	$\gamma_{1}, \ldots, \gamma_{l}$,
	i.e.\ for each $\gamma_{i}$ remove an open
	tubular neighborhood $\nu(\gamma_{i}) \cong \sphere{1} \times \interior \disk{3}$
	and glue
	copies of $\disk{2} \times \sphere{2}$ to the
	new $\sphere{1} \times \sphere{2}$ boundary components
	via the identity map
	$\sphere{1} \times \sphere{2} \rightarrow \sphere{1} \times \sphere{2}$.
	After this surgery we have a compact $4$-manifold $V'$ with $H_{1}(V')=0$,
	and the original disk $\Delta$ survives into $V'$
	in which we will call it $\Delta'$.
	Observe that this ``new'' disk $\Delta'$ is still null-homologous
	in $V'$, since the $3$-manifold is still present after
	the surgery.
    Each circle surgery in a $4$-manifold increases
    the Euler characteristic by 2, thus $\chi(V') = \chi(V) + 2l$.
    By construction, the 4-manifolds $V$ and $V'$
    are cobordant, and so their signatures $\sign(V') = \sign(V)$ agree.
	
	Starting with a knot $K$
	with large enough signature, if there existed
	a null-homologous $\Delta'$,
	since $H_{1}(V') = 0$:
	\begin{equation*}
	    \abs{ \sigma_{K}(-1) + \sign(V') } - \chi(V') + 2
	    =
	    \abs{ \sigma_{K}(-1) + \sign(V) } - (\abs{\chi(V)} + 2l) + 2
	    >
	    0
	\end{equation*}
	which contradicts the inequality in
	\autoref{cor:MT_general_manifold}.
	Therefore $\Delta$ for $K$ cannot exist.
\end{proof}

\begin{remark}
    Earlier sources for results
    in the smooth category
    include Gilmer and Viro's
    \cite{MR603768} version
    of the Murasugi-Tristram inequality for the classical signature
    as stated in
    \cite[Thm.\ 3.1]{yasuhara1996homologyclasses}.
    Our preference for using
    \cite{conway2020stably} in the proof of
    \autoref{thm:no_universal_slicing} comes from
    the result being stated in the topological
    locally flat category.
\end{remark}

\section{Speculation and Questions}
\label{sec:questions}

\subsection{Connection to other conjectures}

An alternative approach to the \SPC is to find a compact 3-manifold $M$
that embeds smoothly in some homotopy 4-sphere $\Sigma^{4}$, but not in $\sphere{4}$.
Notice that if a smooth integral
homology sphere $M$ smoothly embeds in $\Sigma$,
then $M$ is the boundary of a smooth
homology 4-ball \cite[Prop.\ 2.4]{MR3653313}.
However, there is no known example of a 3-manifold $M$ that is the boundary of
a smooth homology 4-ball but that does not embed into $\sphere{4}$.
Both this and the approach in the introduction
are hung up at the homological level.
Further discussion of knots in homology spheres
and concordance in homology cylinders
can be found in, for example, \cite{hom2018knot}, \cite{davis2019concordance}.

\autoref{cor:universal_cover_S4_no_deep_slice} has some relevance to
this which we now discuss (a similar discussion also appears in a comment by
Ian Agol on Danny Calegari's blogpost \cite{scharlemannonschoenflies}).
The unsolved Schoenflies conjecture proposes that if $\mathcal{S} \subset \sphere{4}$
is a smoothly embedded submanifold with $\mathcal{S}$ homeomorphic to $\sphere{3}$,
then $\mathcal{S}$ bounds a submanifold $B \subset \sphere{4}$
that is diffeomorphic to $\disk{4}$.
The \SPC implies the Schoenflies conjecture.  

\begin{question} \label{all embed}
	Does every exotic homotopy 4-ball $\mathcal{B}$ smoothly embed into $\sphere{4}$?
\end{question}

Note that if the answer to \autoref{all embed} is yes,
then the Schoenflies conjecture implies the \SPC and hence the
two conjectures are equivalent: If any homotopy 4-ball
would embed into $\sphere{4}$ and thus,
by the Schoenflies conjecture, would be diffeomorphic to $\disk{4}$,
hence all homotopy 4-balls would be standard,
so all homotopy 4-spheres would be standard.
We have:
\begin{observation}
    If the answer to \autoref{all embed} is yes,
    then no homotopy 4-ball can have deep slice knots.  
\end{observation}

Thus by \autoref{cor:universal_cover_S4_no_deep_slice},
if the answer to \autoref{all embed} is yes,
the approach towards \SPC mentioned in this section would never succeed.
Similarly, there would be no 3-manifold that would smoothly embed into a homotopy 4-sphere
but not into $\sphere{4}$. This is because any such embedding into a homotopy sphere avoids a standard 4-ball and after removing this ball the complement is a homotopy 4-ball which we assume embeds into $S^4$, so this approach to \SPC would also be a dead end.  

\subsection{More questions}

\begin{question}
	Are there any 2-handlebodies $X$ other that
	$\natural^k (\sphere{2} \times \ball{2}), k \ge 0$,
	with the property that all $K \subset \ball{3} \subset X$ that are slice in $X$
	are also slice in $\ball{4}$?
	In other words, are there always deep slice local knots
	when $X \neq \natural^k (\sphere{2} \times \ball{2})$? 
\end{question}

One strategy for answering this question would be to start with a framed link $L$
describing a 2-handlebody other than $\natural^k (\sphere{2} \times \ball{2})$
and to handle-slide $L$ to a new framed link $L'$ that 
contains a knot $K$ that is not slice in $\ball{4}$.
Then this knot $K$ when considered in a 3-ball
$K \subset \ball{3} \subset \partial X$ is an example of such a deep slice knot in $X$.
This strategy fails to find any non-slice knots $K$ (as it must)
for $\natural^k (\sphere{2} \times \disk{2})$ when we start with $L$ being
the 0-framed unlink -- since then all resulting knots $K$
will be ribbon hence slice in $\ball{4}$.  

In view of the $2$-handlebodies constructed in
\autoref{prop:all_slice}, 
one could ask whether this extension of the Norman-Suzuki trick is
the only way to make any knot in the boundary of a manifold
bound an embedded disk:
\begin{question}
	If $X$ is a 2-handlebody with the property that every knot in the boundary
	of $X$ is slice in $X$
	(no assumption on the relative homology class of the disk), 
	does it follow that $X$ decomposes
	as $X = X_0 \# (S^2 \times S^2)$ or $X = X_0 \# (S^2 \widetilde{\times} S^2)$?
	More generally, what about the same question
	without the hypothesis that $X$ be a 2-handlebody?
\end{question}

\bibliographystyle{alpha}
\bibliography{references}

\end{document}